\documentclass{amsart}

\usepackage{amsmath}
\usepackage{amsfonts}
\usepackage{amssymb}

\newtheorem{theorem}{Theorem}[section]

\theoremstyle{definition}
\newtheorem{definition}[theorem]{Definition}
\newtheorem{example}[theorem]{Example}

\theoremstyle{remark}

\numberwithin{equation}{section}

\begin{document}

\title[Quantum Statistical Mechanics and Class Field Theory ]{Quantum Statistical 
Mechanics and Class Field Theory }

%Information for first author
\author{Jorge Plazas}

\address{ Max Planck Institute for Mathematics,
Vivatsgasse 7,
Bonn 53111,
Germany}

\email{plazas@mpim-bonn.mpg.de}

%    General info
%\date{December 20, 2005}


\begin{abstract}
In this short communication we survey some known results relating
noncommutative geometry to the class field theory of number fields.
These results appear within the context of quantum statistical mechanics
where some arithmetic properties of a given number field can be
realized in terms of the structure of equilibrium states of a
quantum statistical mechanical system.

\end{abstract}

\maketitle

\section{Introduction}

In the last ten years some very interesting results relating noncommutative
geometry to class field theory have emerged. The first
instance of this connection was explored by Bost and Connes in
\cite{BostConnes} where they related the class field theory of the
field of rational numbers $\mathbb{Q}$ to the structure of
equilibrium states of a particular quantum statistical mechanical
system. Various results generalizing some aspects of this
construction, based on quantum statistical mechanical systems related
to other number fields, have appeared  since then
(\cite{ALR,Cohen,HL,LR,LF}). In particular the
problem of finding  quantum statistical mechanical systems that
encode the explicit class field theory of quadratic imaginary
extensions of $\mathbb{Q}$  has been solved recently (\cite{CMR1,CMR2}).
The existence of quantum statistical mechanical systems with rich
arithmetical properties opens a new approach to the study of
explicit class field theory using the tools of quantum statistical
mechanics. The purpose of this article is to give a brief introduction to this topic.

I want to thank the organizers of the 2005 summer
school ``Geometric and topological methods for quantum field theory" and
Max Planck Institute for their support.

\section{Basics in class field theory}

The main objects of study in algebraic number theory are number
fields, by definition a number field $K$ is a finite degree
extension of the field of rational numbers $\mathbb{Q}$. Once we fix
an algebraic closure ${\bar K}$ of $K$ we would like to understand
the absolute Galois group $Gal(\bar{K}|K)$, which turns out to be
very difficult even when $K=\mathbb{Q}$. One step in understanding the
group $Gal(\bar{K}|K)$ consists in studying its abelianization
$Gal(\bar{K}|K)^{ab}$, this abelian group is the Galois group of
$K^{ab}$,the maximal abelian extension  of $K$, so we have:
$$
Gal(K^{ab}|K)=Gal(\bar{K}|K)^{ab}
$$

The field $K^{ab}$ may be obtained as the limit over all finite normal
abelian extensions of $K$ in $\bar{K}$. Abelian class field theory
studies these abelian extensions together with the
corresponding Galois groups. One of the main results of the theory characterizes
the group $Gal(K^{ab}|K)$ as a quotient of the group of units of a
topological ring constructed from the field $K$ alone. This
topological ring is the id\`{e}le class group $\mathcal{C}_{K}$.
Let us briefly recall its construction.

Let $K$ be a number field and let $\mathcal{O}_{K}$ be its ring of
integers, i.e. $\mathcal{O}_{K}$ consists of the elements in $K$
which are roots of monic polynomials with coefficients in
$\mathbb{Z}$. There is a one to one correspondence between the prime
ideals of $\mathcal{O}_{K}$ and the finite places of $K$, which are
equivalence classes of nonarchimedean valuations $|\quad |$ on $K$.
For instance every prime number $p\in
\mathcal{O}_{\mathbb{Q}}=\mathbb{Z}$ gives rise to a finite place
corresponding to the $p$-adic absolute value $|\quad |_{p}$ in
$\mathbb{Q}$. Given a finite place $\nu$ of $K$ we denote by
$K_{\nu}$ the completion of $K$ with respect to the metric induced
by $\nu$ and denote by $\mathcal{O}_{\nu}$ the ring of integers of
$K_{\nu}$. For each $\nu$ the ring $\mathcal{O}_{\nu}$ is an open
compact subring of $K_{\nu}$ and we may form the restricted
topological product
$$
\mathbb{A}_{f,K} = \prod_\nu (K_{\nu}:\mathcal{O}_{\nu})
$$
where $\nu$ runs over the finite places of $\mathcal{O}_{K}$. This
restricted product is by definition the subset of the cartesian
product $ \prod_\nu K_{\nu}$ consisting of elements for which all
but a finite number of coordinates lie in the subrings
$\mathcal{O}_{\nu}$. It is given the weakest topology for which the
sets $ \prod_{\nu \in F} K_{\nu} \times  \prod_{\nu \notin F}
\mathcal{O} _{\nu }$, for $F$ a finite set of places, are open. The
topological ring $\mathbb{A}_{f,K}$ is called the ring of finite
ad\`{e}les of $K$. If we add to $\mathbb{A}_{f,K}$ the product of
the completions of $K$ with respect to the infinite places we get
$\mathbb{A}_{K}$, the ring of ad\`{e}les of $K$. Infinite places are
given by equivalence classes of archimedean valuations and
correspond to the $[ K:\mathbb{Q} ]$ different embeddings of $K$ in
$\mathbb{C}$.

Consider now $Gl_{1}(\mathbb{A}_{K})$, the group of units of
$\mathbb{A}_{K}$. $K^{*}$ can be embedded diagonally as a discrete
subgroup of $Gl_{1}(\mathbb{A}_{K})$. The quotient
$\mathcal{C}_{K}=Gl_{1}(\mathbb{A}_{K})/K^{*}$ is called the
id\`{e}le class group of $K$. It turns out that the abelian
extensions of $K$ correspond to the normal subgroups of
$\mathcal{C}_{K}$. In particular the following important
theorem holds (cf. \cite{Gras}, Section II.3.7):

\begin{theorem}
Let $D_K$ be the connected component of the identity in
$\mathcal{C}_{K}$. There is a canonical isomorphism:
$$
\theta: \mathcal{C}_{K} / D_{K} \rightarrow Gal(K^{ab}|K)
$$
\end{theorem}

For the case $K=\mathbb{Q}$ a lot more is known. By the Kronecker Weber
theorem (cf. \cite{Swinnerton-Dyer} Section 20) every abelian extension of 
$\mathbb{Q}$ is contained in a
cyclotomic extension, that is, in some extension of
the form $\mathbb{Q} (\zeta _n)$  where $\zeta _n$ is a
primitive $n$-th root of unity for some $n$. The cyclotomic extension $\mathbb{Q}
(\zeta _n)$ is abelian and its Galois group $Gal(\mathbb{Q} (\zeta
_n)| \mathbb{Q})$ is isomorphic to $(\mathbb{Z} / n \mathbb{Z})^{*}$.
The maximal abelian extension of $\mathbb{Q}$ is the field
$\mathbb{Q}^{cycl}$ obtained by adjoining to $\mathbb{Q}$ all roots
of unity and $\mathcal{C}_{\mathbb{Q}} / D_{\mathbb{Q}} = \hat{ \mathbb{Z}}^{*}$
where $\hat{ \mathbb{Z}} =\underset{\leftarrow n}{\lim} \quad \mathbb{Z} /  n \mathbb{Z}$.

For a general number field $K$, one would also like to know which elements
of $\bar{K}$ generate $K^{ab}$ over $K$ and how the Galois group
$Gal(K^{ab}|K)$ acts on these generators. This is the content of
Hilbert's 12th problem known also as the explicit class field theory
problem. Up to now the only number fields apart from $\mathbb{Q}$
for which a complete solution to the explicit class field theory
problem is known are quadratic imaginary fields, that is fields of the
form $K = \mathbb{Q}(\sqrt{-d})$, where $d \in \mathbb{N}^{+}$ is a positive 
integer. In this case the answer is given by the theory
of complex multiplication which characterizes the field
$\mathbb{Q}(\sqrt{-d})^{ab}$ in terms of data coming from an elliptic
curve (cf. \cite{Serre}). Given an embedding $\mathbb{Q}(\sqrt{-d})
\hookrightarrow \mathbb{C}$ the ring of integers of
$\mathbb{Q}(\sqrt{-d})$ gives rise to a lattice in $\mathbb{C}$ and we
obtain an elliptic $E$ curve as the quotient
of $\mathbb{C}$ by this lattice. The curve $E$ has
a natural abelian group structure and $\mathbb{Q}(\sqrt{-d})^{ab}$
is generated by the values of an  analytic function $F$ on $E$ at
the points of finite order of $E$. The function $F$ is given by  
$\mathcal{P}^{u}$ where
$\mathcal{P}$ is the Weierstrass function of $E$ and $u$ is the
order of the group of automorphisms of $E$.

\section{Basics in quantum statistical mechanics}

Quantum statistical mechanics studies statistical ensembles of
quantum mechanical systems. From a mathematical point of view,
a quantum statistical mechanical system consists of a set of observables $A$
having the structure of a $C^{*}$-algebra and a time evolution of
the system given by a one-parameter group of automorphisms $\sigma_{t}$ of the
algebra of observables, $\sigma_{t} \in Aut(A),\quad t\in \mathbb{R}$.
Classically the state of a system is specified by a probability measure on the phase
space, integration of an observable against this measure gives its mean value 
corresponding to
the particular state. In a quantum
mechanical setting, probability measures are replaced by states on
the algebra of observables. By definition a state on a unital $C^{*}$-algebra $A$ is a
linear map $\varphi : A \rightarrow \mathbb{C}$ satisfying:

\begin{itemize}

\item $\varphi (a^{*}a) \geq 0$ for all  $a \in A$.

\item $\varphi (1)=1$.

\end{itemize}

The appropriate definition of equilibrium states in this context was given by 
Haag, Hugenholtz and Winnink in \cite{Haag}. Given a quantum statistical mechanical 
system $(A,\sigma_t)$ an equilibrium state will 
be a state $\varphi$ on the algebra
$A$ satisfying certain compatibility
condition with respect to the time evolution of the system. 
This condition, known as the $KMS$
condition (after Kubo, Martin and Schwinger)  depends
on a thermodynamic parameter $\beta=\frac{1}{T} $, the inverse temperature of the system. 

\begin{definition} 

Let $(A,\sigma _t)$ be a quantum statistical mechanical system. A state
$\varphi$ on $A$ satisfies the $KMS$ condition at inverse temperature
$0 < \beta < \infty$ if for every $a, b \in A$ there exist a bounded holomorphic function
$F_{a,b}$ on $\{ z\in \mathbb{C} \mid 0< \Im(z) < \beta \}$, continuous
on the closed strip, such that
$$
F_{a,b}(t) = \varphi (a \sigma_t (b)) , \quad F_{a,b}(t + i \beta) = \varphi (\sigma_t (b) a), \quad
\forall t \in
\mathbb{R}
$$

We call such state a $KMS_{\beta}$ state.

\end{definition}

A $KMS_\infty$ state is by definition a weak
limit of $KMS_{\beta}$ states as $\beta \rightarrow \infty$. This definition is stronger 
than the more commonly used definition of ground states as states for which the 
$KMS$ condition holds with $\beta =  \infty$. 
The definition of $KMS_\infty$ states as weak limits of $KMS_{\beta}$ states is better 
behaved and more appropriate for the applications we will describe below.

It turns out (cf. \cite{Bratteli}, Section $5.3$) that for each $\beta$ the set of
$KMS_{\beta}$ states associated to the time evolution $\sigma_t$ is a
compact convex space. We denote by $\mathcal{E} _\beta$ the
space of extremal points of the space of $KMS_ \beta$ states.

\begin{example}

Let $\beta > 0$ and let $A =\mathcal{B}(\mathcal{H})$ be the algebra
of bounded operators on a Hilbert space $\mathcal{H}$. Given a
positive self adjoint operator $H\in A$ we can define a time
evolution on $A$  by
$$
\sigma _t (a) = e^{i t H} a e^{- i t H} ,\qquad a\in A
$$
If the operator $e^{- \beta H}$ is trace class then
$$
\varphi (a) = \frac{1}{Z} Tr (a e^{- \beta H}),\qquad Z = Tr (e^{- \beta H})
$$
is a $KMS_{\beta}$ state.
\end{example}

More generally, given a quantum statistical mechanical system $(A,\sigma _t)$  we can look for
representations of the algebra $A$ as an algebra of bounded operators in a Hilbert
space $\mathcal{H}$. Then, given a positive self adjoint operator $H$ such that the time 
evolution
 $\sigma_{t}$ has the form above, we call
$$
Z(\beta) = Tr (e^{- \beta H})
$$
the partition function of the system. The poles of this function
correspond to the critical temperatures at which phase
transitions of the system may take place.

Consider a quantum statistical mechanical system $(A,\sigma _t)$. A group $G \subset Aut(A)$
such that
$$
\sigma_t g = g \sigma_t \qquad \forall g \in G, \quad \forall t \in
\mathbb{R}
$$
is called a symmetry group of the system. $G$ then acts on the space of $KMS_ \beta$ states
for any $\beta$ and hence on $\mathcal{E}_{\beta}$. Inner automorphisms
coming from unitaries invariant under the time evolution act trivially on equilibrium states.

For some of the applications we will describe below 
it is also important to consider symmetries of the system induced by
endomorphisms of the algebra $A$ (\cite{ConnesMarcolli, CMR1,CMR2}). An endomorphism of $A$ 
is by definition a $*$-morphism $\rho:A\rightarrow A$ commuting with $\sigma_{t}$. Such an
endomorphism acts on the set of $KMS_{\beta}$ states $\varphi$ for which $\varphi(\rho(1))
\neq 0$ by $\rho^{*}\varphi = \frac{1}{\varphi (\rho(1))}\varphi \circ \rho$ .
Inner endomorphisms coming from isometries invariant under the time evolution act
trivially on equilibrium states. Symmetries induced by endomorphisms were introduced 
in the context of super selection sectors developed by Doplicher, Haag and roberts.

\section{The Bost-Connes system}

In \cite{BostConnes}  Bost and Connes constructed a remarkable quantum
statistical mechanical system $(A, \sigma_t)$ in which the structure
of equilibrium states is related in a deep way with the class field
theory of $\mathbb{Q}$. The group 
$\mathcal{C}_{\mathbb{Q}} / D_{\mathbb{Q}} = \hat{ \mathbb{Z}}^{*}$
acts as symmetries of this system and the algebraic numbers generating the maximal abelian
extension of $\mathbb{Q}$ can be recovered as values of the
equilibrium states at zero temperature on the observables
corresponding to an arithmetic subalgebra of $A$. What is more
remarkable about this system is the fact that  the action of the
Galois group on this values commutes with the action of
$\mathcal{C}_{\mathbb{Q}} / D_{\mathbb{Q}} = \hat{ \mathbb{Z}}^{*}$ on the equilibrium states.
In this section we describe this system and its main features.

Let $A$ be the $C^{*}$-algebra generated by two sets of elements
$\{ e(r) \mid r\in \mathbb{Q}/ \mathbb{Z} \}$ and $\{ \mu_n \mid n \in
\mathbb{N}^{+} \}$  with relations:

\begin{enumerate}
\item $\mu_{n}^{*} \mu_n = 1 \qquad \forall n \in \mathbb{N}^{+}$

\item $\mu_{n}\mu_{k}=\mu_{nk} \qquad\forall n, k \in \mathbb{N}^{+}$
\item $e(0)= 1 ,  \quad e(r)^{*}=e(-r), \quad  e(r)e(s)=e(r+s) \qquad \forall
r, s \in \mathbb{Q}/\mathbb{Z}$

\item $\mu_{n} e(r) \mu_{n}^{*} = \frac{1}{n} \sum_{ns=r} e(s) \qquad  \forall n
\in \mathbb{N}^{+},  \forall r \in \mathbb{Q}/\mathbb{Z}$

\end{enumerate}
Define a time evolution on the $C^{*}$-algebra $A$ by taking
$$
\begin{array}{cc}
\sigma_t (\mu_{n}) = n^{i t} \mu_{n},\qquad  n \in \mathbb{N}^{+}, \quad t\in \mathbb{R} \\
\quad \sigma_t (e(r)) = e (r), \qquad \qquad   r \in \mathbb{Q}/\mathbb{Z},  \quad t\in 
\mathbb{R} \\
 \end{array}
$$
We will refer to the quantum statistical mechanical system $(A,\sigma_t)$ as the Bost-Connes 
system. We summarize the main results of \cite{BostConnes} about the structure of this 
system. As above we denote by $\mathcal{E}_{\infty}$ the set of extremal $KMS_{\infty}$
states.

\begin{theorem}[Bost, Connes]\quad

\begin{itemize}

\item The group $\mathcal{C}_{\mathbb{Q}} / D_{\mathbb{Q}} = \hat{
\mathbb{Z}}^{*}$ acts as a group of symmetries of the system $(A,
\sigma_t)$.

\item Let $\mathcal{A}_{\mathbb{Q}}$ be the $\mathbb{Q}$-subalgebra
of $A$ generated over $\mathbb{Q}$ by the sets of elements
$\{ e(r) \mid r\in \mathbb{Q}/ \mathbb{Z} \}$ and $\{ \mu_n, \mu_{n}^{*} \mid n \in
\mathbb{N}^{*} \}$.
Then for every $\varphi \in \mathcal{E}_{\infty}$ and every $ a\in
\mathcal{A}_{\mathbb{Q}}$ the value $\varphi (a)$ is algebraic over
$\mathbb{Q}$. Moreover for any $\varphi \in \mathcal{E}_{\infty}$ one has
$\varphi ( \mathcal{A}_{\mathbb{Q}}) \subset \mathbb{Q}^{ab}$ and $ \mathbb{Q}^{ab}$ is
generated by numbers of the  form $\varphi (a)$ with $\varphi \in \mathcal{E}_{\infty}$ and
$ a\in \mathcal{A}_{\mathbb{Q}}$.

\item For all $\varphi \in \mathcal{E}_{\infty}$, $\gamma \in
Gal(\mathbb{Q}^{ab}|\mathbb{Q})$ and $a\in \mathcal{A}_{\mathbb{Q}}$
one has
$$
\gamma \varphi (a) = \varphi (\theta ^{-1} (\gamma) a)
$$
where $\theta: \hat{\mathbb{Z}}^{*} \rightarrow
Gal(\mathbb{Q}^{ab}| \mathbb{Q}) $ is the class field theory
isomorphism.
\end{itemize}

\end{theorem}

The system $(A,\sigma_t)$ was originally introduced by Bost and Connes in the context
of Hecke algebras. The inclusion of rings
$\mathbb{Z} \subset \mathbb{Q}$ induces an inclusion matrix groups $P^{+}_{\mathbb{Z}}\subset
P^{+}_{\mathbb{Q}}$ where $P^{+}_{\mathbb{Q}} =\{ \left(
                              \begin{array}{cc}
                                 1 & b  \\
                                 0 & a \\
                              \end{array}
                            \right) | b\in \mathbb{Q}, a \in \mathbb{Q} ^{*}_{+} \}$
and $P^{+}_{\mathbb{Z}} =\{ \left(
                              \begin{array}{cc}
                                 1 & c  \\
                                 0 & 1 \\
                              \end{array}
                            \right) | c\in \mathbb{Z} \}$. This inclusion of matrix groups
is almost normal in the sense that the orbits of $P^{+}_{\mathbb{Z}}$ acting on the left on 
the set of right cosets $P^{+}_{\mathbb{Q}} / P^{+}_{\mathbb{Z}}$ are finite and viceversa. 
One can then form a convolution algebra $\mathcal{H}(P^{+}_{\mathbb{Z}} , P^{+}_{\mathbb{Q}})$ 
given by finitely supported $\mathbb{C}$-valued functions on
$ P^{+}_{\mathbb{Z}} \backslash P^{+}_{\mathbb{Q}} / P^{+}_{\mathbb{Z}}$, this algebra is 
called the Hecke algebra of the pair ($P^{+}_{\mathbb{Z}},P^{+}_{\mathbb{Q}})$.
$\mathcal{H}(P^{+}_{\mathbb{Q}}; P^{+}_{\mathbb{Z}} )$ can be completed
to a $C^{*}$-algebra with a natural time evolution. It is shown in \cite{BostConnes} that in
this way one obtains the system $(A,\sigma_t)$ described above.

Let $(\epsilon_{n})_{n\in \mathbb{N}^{+}}$  be
the canonical basis of the Hilbert space $l^{2}(\mathbb{N}^{+})$.
For any element $\gamma \in Gal(\mathbb{Q}^{ab}|\mathbb{Q})$  one can define a
representation $\pi_{\gamma}$ of the algebra $A$ in $l^{2}(\mathbb{N}^{+})$ by
$$
\begin{array}{cc}
\pi_{\gamma}(\mu_{n})\epsilon_{k} = \epsilon_{n k},\qquad  n,k \in \mathbb{N}^{+} \\
\pi_{\gamma}(e(r)) \epsilon_{k} = \gamma (e^{2 \pi i k r})\epsilon_{k}, \qquad  k \in
\mathbb{N}^{+}, r \in \mathbb{Q}/\mathbb{Z} \\
 \end{array}
$$
For any of these representations the time evolution of the system is
implemented by the operator
$ H \epsilon_{k} = (\log k) \epsilon_{k}, \quad  k \in \mathbb{N}^{+} $
and so the partition function of the system is the Riemann zeta function
$$
Z(\beta) = Tr(e^{- \beta H}) =
\sum_{k} \frac{1}{k^{\beta}} = \zeta(\beta).
$$

The Bost-Connes system exhibits a phenomenon called spontaneous symmetry breaking.
 For $0 < \beta \leq 1$  the system $(A,\sigma_t)$ admits only one
equilibrium state and the symmetry group of the system acts
trivially on this state. For $\beta > 1 $ the space
$\mathcal{E}_{\beta}$ is parameterized by $\hat{\mathbb{Z}}^{*}$ and
the symmetry group acts transitively on this space.

\section{Some generalizations to other number fields}

The work of Bost and Connes has inspired several authors to construct quantum
statistical mechanical systems associated to other number fields. These
constructions generalize in different directions some of the results in \cite{BostConnes}.

In \cite{LR} M. Laca and I. Raeburn realized the algebra of
observables of the Bost-Connes system as a semigroup crossed product
algebra. Given a discrete group $\Gamma$ its group algebra
$\mathbb{C}\Gamma$ can be completed to a $C^{*}$-algebra
$C^{*}(\Gamma)$ by considering all unitary irreducible
representations of $\mathbb{C}\Gamma$ as an algebra of operators on
some Hilbert space. If a semigroup $S$ acts by endomorphisms on the
algebra $C^{*}(\Gamma)$ one can twist the product and the
convolution on $C^{*}(\Gamma)$ by the action of $S$ getting a
crossed product $C^{*}$-algebra $C^{*}(\Gamma) \rtimes S$.  The
$C^{*}$-algebra of the Bost-Connes system is then given by a
semigroup crossed product $C^{*}(\mathbb{Q} / \mathbb{Z}) \rtimes
\mathbb{N}^{+}$ where the action of $ \mathbb{N}^{+}$ on
$C^{*}(\mathbb{Q} / \mathbb{Z})$ is a right inverse to the action
coming from the natural multiplication $ (n,r) \mapsto nr$ where
$n\in \mathbb{N}^{+}$ and $r\in \mathbb{Q} / \mathbb{Z}$.

Arledge, Laca and Raeburn considered in \cite{ALR} crossed
products of the form $C^{*}(K / \mathcal{O}) \rtimes
\mathcal{O}^{\times}$ for $K$ an arbitrary number field and
$\mathcal{O}$  its ring of integers. They characterize the faithful
representations of these algebras and realized them as Hecke
algebras.

As in the case of $\mathbb{Q}$ one has that for any number field $K$
the ring inclusion $\mathcal{O} \subset K$ induces an almost normal
inclusion of matrix groups $P^{+}_{\mathbb{\mathcal{O}}}\subset
P^{+}_{K}$. Harari and Leichtnam studied in \cite{HL} a
quantum statistical mechanical system corresponding to the Hecke
algebra $\mathcal{H}(P^{+}_{K}; P^{+}_{\mathcal{O}})$. The
equilibrium states of this system share many of the properties of
the Bost-Connes system. The group $\hat{\mathcal{O}}^{*}= \prod_\nu
\mathcal{O}_{\nu}^{*}$, where $\nu$ runs over the finite places of
$K$, acts as a symmetry group of the system and the partition function
of the system is the Dedekind zeta function of the number field $K$
$$
\zeta_K(s) = \sum_{\mathfrak{a}} \frac{1}{N( \mathfrak{a})^{s}} =  
\prod_{\mathfrak{p}}\frac{1}{1 -N(\mathfrak{p})^{-s}}
$$
with a finite number of factors removed from the product. In this
expression $\mathfrak{a}$ runs over the integral ideals of
$\mathcal{O}$, $\mathfrak{p}$ runs over the prime ideals of
$\mathcal{O}$ and   $N( \mathfrak{a})=[\mathcal{O} : \mathfrak{a}]$
is the absolute norm of the ideal $\mathfrak{a}$. This system
exhibits spontaneous symmetry breaking at the pole of this function;
for $0 < \beta \leq 1$  the system admits only one equilibrium state
and the symmetry group of the system acts trivially on this state.
For $ \beta > 1$ the space $\mathcal{E}_{\beta}$ is parameterized by
$\hat{\mathcal{O}}^{*}$ and the symmetry group acts freely and
transitively on this space. The results of \cite{HL} hold for
arbitrary global fields. These include number fields and finite
separable extensions of the field of rational functions over a
finite field.

A quantum statistical mechanical system similar to the one of \cite{HL} was introduced
in \cite{Cohen} by Cohen. Her system has the advantage of recovering the full Dedekind zeta
function of the number field $K$ as partition function of the system.

Laca and van Frankenhuijsen studied in \cite{LF} the Hecke algebra corresponding to the 
almost normal inclusion of matrix groups
$$
\left(
                              \begin{array}{cc}
                                 1 & \mathcal{O}  \\
                                 0 & \mathcal{O}^{*} \\
                              \end{array}
                            \right)
                \subset
                \left(
                              \begin{array}{cc}
                                 1 & K \\
                                 0 & K^{*} \\
                              \end{array}
                            \right)
$$
for an arbitrary number field $K$ with ring of integers
$\mathcal{O}$. The equilibrium states of the corresponding system
are analyzed in the case of fields with class number one. For
quadratic imaginary fields of class number one the symmetry group of
the system is isomorphic to the Galois group of the maximal abelian
extension of the field and the action of this group on extremal
equilibrium states is transitive.

\section{Fabulous states for number fields}

The work of Bost and Connes together with its various
generalizations opens the possibility of approaching the problem of explicit
class field theory within the framework of noncommutative
geometry. Following these lines one would like to have for an
arbitrary number field $K$ a quantum statistical mechanical system which fully
incorporates its class field theory, the notion of fabulous states for number fields,
introduced in \cite{ConnesMarcolli} by Connes and  Marcolli, encodes this aim.

Given a number field $K$ together with an embedding $K\hookrightarrow \mathbb{C}$
the ``problem of fabulous states" asks for the construction of a quantum statistical
mechanical system $(A,\sigma_t)$ such that:

\begin{itemize}

\item The group $\mathcal{C}_{K} / D_{K}$ acts as a group of symmetries of the system $(A,
\sigma_t)$.

\item There exists a $K$-subalgebra
$\mathcal{A}_{K}$ of $A$ such that for every
$\varphi \in \mathcal{E}_{\infty}$ and every $ a\in
\mathcal{A}_{K}$ the value $\varphi (a)$ is algebraic over
$K$. Moreover $K^{ab}$ is generated over $K$ by numbers of
this form.

\item For all $\varphi \in \mathcal{E}_{\infty}$, $\gamma \in
Gal(K^{ab}|K)$ and $a\in \mathcal{A}_{K}$
one has
$$
\gamma \varphi (a) = \varphi (\theta^{-1}(\gamma) a)
$$
where $\theta: \mathcal{C}_{K} / D_{K} \rightarrow
Gal(K^{ab} |K)$ is the class field theory
isomorphism.
\end{itemize}

The ground states of such a system are
refereed as ``fabulous" in view of its rich arithmetical properties.

\section{$\mathbb{Q}$-lattices, quadratic extensions and complex multiplication}

For the case of quadratic imaginary fields $K = \mathbb{Q}(\sqrt{-d})$, $d \in \mathbb{N}^{+}$,
the problem of finding a quantum statistical mechanical system
having all the properties discussed in the last section has been recently solved by Connes,
Marcolli and Ramachandran \cite{CMR1,CMR2}. Their results rely deeply on the arithmetic
properties of $\mathbb{Q}$-lattices which were studied in \cite{ConnesMarcolli}.

Given a number field $K$ with $[K,\mathbb{Q}]=n$ there is an
embedding $K^{*} \hookrightarrow Gl_{n}(\mathbb{Q})$, this induces
an embedding $Gl_1(\mathbb{A}_{K}) \hookrightarrow
Gl_{n}(\mathbb{A}_{\mathbb{Q}})$. It is therefore important to study
$Gl_n$ analogs of the Bost-Connes system. In \cite{ConnesMarcolli}
Connes and Marcolli constructed a quantum statistical mechanical
system with group of symmetries $Gl_{2}(\mathbb{A}_{\mathbb{Q}, f})$
and analyzed in detail the arithmetical properties of its
equilibrium states. The central notions introduced Connes and
Marcolli in \cite{ConnesMarcolli} are those of
$\mathbb{Q}$-lattices and commensurability.

A $n$ dimensional $\mathbb{Q}$-lattice is by definition given by a pair $(\Lambda, \phi)$ 
where $\Lambda$ is a lattice in $\mathbb{R}^{n}$ and  $\phi: \mathbb{Q}^{n} / \mathbb{Z}^{n} 
\rightarrow \mathbb{Q} \Lambda  / \Lambda$
is a homomorphism of abelian groups. Two $n$ dimensional
$\mathbb{Q}$-lattices $(\Lambda_1, \phi_1)$ and $(\Lambda_2, \phi_2)$
are commensurable if $\mathbb{Q}\Lambda_1 = \mathbb{Q} \Lambda_2$
and $ \phi_1=  \phi_2$ modulo $\Lambda_1 + \Lambda_2$. The relation
of commensurability is an equivalence relation on the set of $n$
dimensional $\mathbb{Q}$-lattices. The space of commensurability classes of
$n$-dimensional $\mathbb{Q}$-lattices is described in terms of a noncommutative 
$C^{*}$-algebra $C^{*}(\mathcal{L}_{n})$.

A one dimensional $\mathbb{Q}$-lattice can be rescaled by
multiplying it by a positive scale factor so the multiplicative
group $\mathbb{R}^{*}_{+}$ acts on the set of one dimensional
$\mathbb{Q}$-lattices. The $C^{*}$-algebra
$C^{*}(\mathcal{L}_{n}/\mathbb{R}^{*}_{+})$ corresponding to the
space of commensurability classes of one dimensional
$\mathbb{Q}$-lattices up to scaling has a natural time evolution and
a natural choice of an arithmetic structure given by a
$\mathbb{Q}$-subalgebra of
$C^{*}(\mathcal{L}_{n}/\mathbb{R}^{*}_{+})$. In this way one
recovers the Bost-Connes system \cite{ConnesMarcolli}.

The group  $\mathbb{C}^{*}$ acts by rescaling  on the set of two dimensional
$\mathbb{Q}$-lattices. As in the one dimensional case the $C^{*}$-algebra
$C^{*}(\mathcal{L}_{n}/\mathbb{C}^{*})$ corresponding to the space
of commensurability classes of two dimensional $\mathbb{Q}$-lattices
up to scaling has a natural time evolution and
a natural choice of an arithmetic structure. The structure of equilibrium states of
the corresponding quantum statistical mechanical system
is related to the Galois theory of the field of modular functions. We refer to
\cite{ConnesMarcolliSurvey} and \cite{Marcolli} for a survey of these results.

Consider now a quadratic imaginary field $K = \mathbb{Q}(\sqrt{-d})$, $d \in \mathbb{N}^{+}$,
with ring of integers $\mathcal{O}$. A one
dimensional $K$-lattice is given by a pair $(\Lambda, \phi)$ where $\Lambda$ is
 a finitely generated $\mathcal{O}$-submodule of $\mathbb{C}$ with 
 $\Lambda \otimes_{\mathcal{O}} K \cong K$
and  $\phi: K / \mathcal{O} \rightarrow K \Lambda  / \Lambda$ is a
homomorphism of abelian groups. Two one dimensional $K$-lattices
$(\Lambda_1, \phi_1)$ and $(\Lambda_2, \phi_2)$ are commensurable if
$K \Lambda_1 = K \Lambda_2$ and $ \phi_1=  \phi_2$ mod $\Lambda_1 +
\Lambda_2$,  the relation of commensurability is  an equivalence
relation on the set of one dimensional $K$-lattices. A one
dimensional $K$-lattice is in particular a two dimensional
$\mathbb{Q}$-lattice and two one dimensional $K$-lattices are
commeasurable if and only if the underlying $\mathbb{Q}$-lattices
are commeasurable. The space of commensurability classes of one
dimensional $K$-lattices gives rise to a quantum statistical
mechanical system whose ground states are fabulous states for the
field $K$  \cite{CMR1,CMR2}. The explicit class field theory of
the quadratic imaginary fields can thus be fully recovered in the form
discussed above.

The spaces associated to commensurability
classes of $\mathbb{Q}$-lattices can be interpreted as
noncommutative versions of Shimura varieties. Classical Shimura
varieties have a rich structure coming from a natural action
of ad\`{e}lic groups on them, they appear in the context of \cite{ConnesMarcolli, CMR1, CMR2} as spaces 
of extremal  $KMS$ states at low temperature. In \cite{Ha} Ha and Paugam developed further this point of view 
constructing quantum statistical mechanical systems associated to arbitrary Shimura
varieties. Given a general number field $K$ the results in \cite{Ha} lead in particular to the existence 
of a quantum statistical mechanical system whose group of symmetries is $\mathcal{C}_{K} / D_{K}$ and whose 
partition function is the Dedekind zeta function $\zeta_K(s)$.

The first case for which there is not yet a complete solution to the explicit class field theory 
problem is the case of real quadratic fields, $K = \mathbb{Q}(\sqrt{d})$, where 
$d \in \mathbb{N}^{+}$ is a square free positive integer. In \cite{manin}  Manin proposed the use 
of noncommutative geometry as a geometric framework for the study of abelian class field 
theory of real  quadratic fields. This is the so called ``real multiplication program". 
The main idea is that noncommutative tori may play a role 
in the study of real quadratic fields analogous to the role played 
by elliptic curves in the study of imaginary quadratic fields. 
Some possible relations between Manin's program and $\mathbb{Q}$-lattices are discussed 
in \cite{Marcolli}.

\bibliographystyle{amsalpha}

\end{document}